\documentclass[10pt,draft,reqno]{amsart}
      \makeatletter
     \def\section{\@startsection{section}{1}%
     \z@{.7\linespacing\@plus\linespacing}{.5\linespacing}%
     {\bfseries
     \centering
     }}
     \def\@secnumfont{\bfseries}
     \makeatother
\newtheorem{theorem}{Theorem}[section]
\newtheorem{lemma}[theorem]{Lemma}

\newtheorem{corollary}[theorem]{Corollary}
\theoremstyle{definition}

\theoremstyle{remark}
\newtheorem{remark}[theorem]{Remark}
\numberwithin{equation}{section}
\def \a{{\alpha}}
\def \b{{\beta}}
\def \D{{\Delta}}

\def \d{{\delta}}
\def \e{{\varepsilon}}

\def \g{{\gamma}}

\def \k{{\kappa}}
\def \l{{\lambda}}

\def \o{{\omega}}

\def \E{{\bf E}\, }

\def \N{{\bf N}}

\def \P{{\bf P}}

\def \qq{{\qquad}}
\def \R{{\bf R}}

\def \dd{{\rm }}

\def \noi{{\noindent}}

\def\E{{\mathbb E}}

\def\P{{\mathbb P}}

\def\R{{\mathbb R}}

\def\Z{{\mathbb Z}}

\def\N{{\mathbb N}}

at  8 pt
\scrollmode

\font\gsec= cmb10 at 10 pt
   at 10,8  pt
  \scrollmode
\font\gsec= cmb10 at 9  pt

 \title[Instants of small amplitude of  Brownian motion]{Instants of small amplitude of  Brownian motion  
  and application to    the 
 Kubilius model} 
 
  \author{ Michel Weber 
}
\address{ Michel Weber: IRMA, Universit\'e
Louis-Pasteur et C.N.R.S.,   7  rue Ren\'e Descartes, 67084
Strasbourg Cedex, France. }
\email{michel.weber@math.unistra.fr
}
 \urladdr{http://www-irma.u-strasbg.fr/$\sim$weber/}
\begin{document}
 \maketitle
   \begin{abstract}Let $W(t), t\ge 0$ be standard Brownian motion.  We study the size of the time intervals  $I$  which are admissible for
the long range  of slow increase, namely given a real $z>0$, 
 $$ \sup_{t\in I}{|W(t)|\over \sqrt t} \le z, $$ 
and we estimate their number of occurences. We obtain optimal results in terms of class test functions  and,  by means of the
quantitative Borel-Cantelli
  lemma, a fine frequency result concerning their occurences. Using Sakhanenko's invariance principe to transfer  the
results to the   Kubilius model, we   derive applications  to the
prime number divisor function. We obtain refinements of some results recently proved by Ford and Tenenbaum in \cite{FT}. 
\end{abstract}
\maketitle
 \vskip 0,2cm  \noi {\gsec 2010 AMS Mathematical Subject Classification}:  Primary 60J65, 60G17 ; Secondary 60G15.
    \par\noi  
{{\gsec Keywords}: Brownian motion, Ornstein-Uhlenbeck process, small deviations, quantitative Borel-Cantelli
  lemma, Sakhanenko's invariance principe, frequency results, class test functions, Kubilius model, prime  divisor function.

 \section{ Introduction-Main results} Let $W(t), t\ge 0$ be standard Brownian motion.   Let $z$ be some positive real.     The study
of the number of occurences of the
 time intervals  $I$ for which 
 $$ \sup_{t\in I}{|W(t)|\over \sqrt t} \le z, $$ 
is the first motivation of  this work. In a second step, we will derive applications for the Kubilius model in number theory.   More
precisely, 
 let $f:[1,\infty)\to\R^+$ be here and throughout   a non-decreasing function such that $f(t)\uparrow \infty$   with
$t$  and
\begin{equation}\label{f}
  f(t)= o_\rho(t^\rho). 
\end{equation}   We will  consider intervals of type 
 $I= [  N ,   N  f(N )]$. We   essentially examine the case $N=e^k$, $k=1,2,\ldots$. The study made can be extended 
with no difficulty to more general  
geometrically increasing sequences, but this aspect will be not developed. Put 
 \begin{equation} \label{akw} A_k(f,z)=\Big\{\sup_{e^k\le t\le e^kf(e^k) } {|W(t)|\over \sqrt t}<z\Big\} , \qq k=1,2,\ldots
\end{equation}
  Let $  U(t)= W(e^t)e^{-t/2}, t\in \R $ be the
Ornstein-Uhlenbeck process. It will be more convenient to work with $U$ instead of $W$. Observe
that$$  A_k(f,z)  = \Big\{\sup_{k\le s\le k +\log f(e^k) }|U(s)|\le
z\Big\}.    $$
And so  as $U$ is   stationary 
$$  \P\{A_k(f,z) \} = \P\Big\{\sup_{0\le s\le \log f(e^k) }|U(s)|\le
z\Big\}.    $$   
   
  We say that
$f\in
\mathcal U_z$ whenever
 $\P\big\{ \limsup_{k\to \infty}A_k(f,z)  \big\}=0 $, and that $f\in \mathcal V_z$ if
 $\P\big\{ \limsup_{k\to \infty}A_k(f,z)  \big\}=1 $.  By the   $0$-$1$ law (since $U$ is strongly mixing), the latter probabilities
can only be 
$0$ or
$1$.   
  
Notice that if
$f\in
\mathcal U_z$, then with probability one 
$$J(f ): = \liminf_{k\to \infty}
\sup_{k\le s\le k +\log f(e^k)}|U(s)|>z, $$
whereas  $J(f)\le z$, almost surely if $f\in \mathcal V_z$. In the latter case, it makes sense to
   estimate the   size of the   counting function
$$N_n(f,z)= \sum_{k=1}^n\chi_{A_k(f,z)}\qq n=1,2,\ldots  $$ 
Naturally  this has to be done with respect to the corresponding  means $   \nu _n(f,z):= \E N_n(f,z)$.

\smallskip\par  We shall
first   characterize  the classes $\mathcal U_z$ and $\mathcal V_z$ by means of a simple convergence criterion, and complete our
characterization  by  including a frequency  result concerning the class $\mathcal V_z$.    
\begin{theorem} \label{fr1}   There exists   $\l(z)>0$ with $\l(z)\sim {\pi^2\over 4z^2} $ as $z\to 0 $, such that if $\Sigma(f) = \sum_k
f(e^k)^{- \l(z)}$, then
 $$f\in \mathcal U_z\quad (\hbox{resp. $\in \mathcal V_z$})\quad \Longleftrightarrow \quad \Sigma(f) <\infty\quad
 (\hbox{resp. $=\infty$}).$$
 Further   for any $a>3/2$,  
$$N_n(f ,z) \buildrel{a.s.}\over{=} \nu _n(f ,z)+ \mathcal O \big(  \nu^{1/2} _n(f ,z)  
\log^a \nu _n(f ,z) \big).$$
   And  there are positive constants $K_1(z),K_2(z)$depending on $z$ only, such that for all $n$  
$$K_1(z)  \le  { \nu _n(f ,z)\over \sum_{k=1}^n  f(e^k)^{- \l(z)}} \le  K_2(z).$$
  \end{theorem}
The critical   value  $\l(z)$ is   the smallest eigenvalue in the Sturm-Liouville
equation (\ref{sl/e}).     See section 2.
\smallskip\par
  The class of functions $f_c(t)=  \log^c t $, $c>0$, is of  special interest  in view of applications to the
Kubilius model.   We deduce from Theorem \ref{fr1}:
\begin{corollary} \label{J}    If $  c  > 1/\l(z) $, then  $f_c\in \mathcal U_z$ whereas  $f_c\in \mathcal V_z$ if $0<c  \le 1/\l(z) $.
Further, for any
$0<c  \le 1/\l(z) $ and
$a>3/2$, 
 \begin{eqnarray*} N_n(f_c,z)&\buildrel{a.s.}\over{=}&\nu _n(f_c,z)+ \mathcal O \Big(  \nu^{1/2} _n(f_c,z)  
\log^a \nu _n(f_c,z) \Big) 
 . 
\end{eqnarray*}
 And for all $n$  
$$K_1(z)  \le  { \nu _n(f_c,z)\over \sum_{k=1}^n k^{-c\l(z)}}  \le  K_2(z).$$
  \end{corollary}

    Accordingly, if 
\begin{equation}\label{Iw} I(f)  :=\liminf_{k \to \infty} \sup_{e^k\le t\le e^kf(e^k) } {|W(t)|\over
\sqrt t}, 
\end{equation}
  then $\P\{  I(f_c) \le z \}=1$ if and only if $0<c  \le 1/\l(z) $.  This is clear  in view of (\ref{akw}).
Noticing  that $I(f)\le I(g)$ whenever $f(N)\le g(N)$ for all $N$ large,  we therefore also deduce 
\begin{corollary}\label{W} We have 
  $\P\{   I(f_c)\le z \}=1$  if and only if    $0<c  \le 1/\l(z) $.
And $\P\{  I(f) =\infty \}=1$  if $f (t)\gg_c   f_c(t)$ for all $c$.
\end{corollary}
\begin{remark}   This slightly improves upon Theorem 3 in  \cite{FT}, where it was shown that 
$\P\{  I(f) <\infty \}=1$ if $f(N)= (\log N)^b$ for some $b>0$, whereas   $\P\{  I(f) =\infty \}=1$, if $f(N)= (\log N)^{b(N)}$ with
$b(N)\to \infty$ with $N$.
 \end{remark}
 
In \cite{FT}, the behavior   of corresponding functionals $I_f$ for sums of independent random
variables (assuming only second absolute moments) was also considered.
 In this direction, we will also establish the following result for sums of independent random variables.
\begin{theorem} \label{ind} Let $\{X_j, j\ge
1\}$ be independent centered random variables.    Assume that $ $   for some $\a > 2$,  
\begin{equation}\label{alpha} \sum_{j\ge 1}\E |X_j|^\a=\infty \qq {\it and}\qq v= \sup_{j\ge 1} {\E |X_j|^\a\over \E |X_j|^2}<\infty.
\end{equation}
Let   
$Z_n= X_1+\ldots+X_n$, $z_n^2=\E Z_n^2$,  $J_n=\big\{ j: n\le z_j^2\le nf(n)\big\}$. Then    there exists a Brownian motion $W$ such that 
$$ \liminf_{k\to \infty}\sup_{e^k\le z_j^2\le e^kf(e^k)}{|Z_j|\over z_j}\buildrel{a.s.}\over{=}\liminf_{k\to \infty}\sup_{e^k\le
z_j^2\le e^kf(e^k)}{|W(z_j^2)|\over s_j}. 
$$
 almost surely. 
  In particular  if $c<1/\l(z)$, then
$$\liminf_{k\to \infty}\sup_{e^k\le z_j^2\le e^kf_c(e^k)}{|Z_j|\over z_j}\le z, \qq {\it almost\ surely}. $$
 \end{theorem}
Notice   in the iid case  that assumption (\ref{alpha}) simply reduces to the integrability condition $\E |X_1|^\a<\infty$ for some
$\a>2$. 
\bigskip\par Now introduce the truncated prime divisor function
  $\o(m,t)=\#\{ p\le t : p|m\}$. Here and throughout we reserve the letter $p$ to denote some arbitrary prime number. Put 
\begin{equation}\label{rho}\rho (m,t) :={ | \o(m,t)- \log\log t |\over \sqrt{\log\log  t}}. 
\end{equation}
 The local variations of  $\rho (m,t) $
were recently investigated by Ford and Tenenbaum, who obtained in   \cite{FT}, after a careful study of the size of intervals of slow growth for general sums of independent random
variables,   quite elaborated asymptotic estimates, on the basis of  the 
approximation formula (\ref{km}).  
 The  results   concern the functional
\begin{equation} \label{funct}   \max_{N\le \log\log  t\le Nf(N)} \rho(m,t).
\end{equation}
 Let  $f(m)$, $g(m)$ be  increasing and   tending to infinity with $m$. 
 It is notably proved
(\cite{FT}, Theorem 5) that
  if  
$g(m)\le (\log\log  m)^{1/10}$  and $f(N)= (\log N)^{\xi(N)}$ where $\xi(N)\to \infty$ sufficiently slowly so that
$f(N)\le N$, then 
\begin{equation} \label{ft}  \min_{ g(m)\le N\atop Nf(N)\le \log\log  m} \max_{N\le \log\log  t\le Nf(N)} \rho(m,t)\to \infty
\end{equation}
along a set of integers $m$ of natural density $1$. 
\smallskip\par 
Further, if $f(N)= (\log N)^c$ and $g^2(m)(\log g^2(m))^c\le \log\log m$ for $m$ large, then on a set of integers $m$ of density $1$, we
have
\begin{equation} \label{ft1}   \min_{ g(m)\le N \le g^2(m)} \max_{N\le \log\log  t\le Nf(N)} \rho(m,t)\le 30 \sqrt{1+c}.
\end{equation}
This provides informations on the size of intervals  which are admissible for the long range  of slow increase, in terms of the
natural density on the integers. For instance $g(m)=\sqrt{(\log\log m)  /(\log\log\log m)^c}$ is suitable. The principle followed in the 
proofs consists with modifying the proofs of the preliminary results on the size of intervals of slow growth for   sums of
independent random variables for the particular sequence $\{T_n, n\ge 1\}$ (section 2) and next to apply approximation
formula (\ref{km}).
\bigskip\par
Here we will proceed slightly differently. As we have optimal results on instants of small amplitude of   Brownian motion, we
directly compare the functionals (\ref{funct}) with analogous functionals of   Brownian motion by means of Sakhanenko's invariance
principle (Lemma \ref{sha1}). This is done in Theorem \ref{outom} below. This allows to transfer our previous results to truncated prime
divisor function, not fully naturally, but sufficiently much to get
  new quite sharp results. More precisely, let
$0<M_1(x)<M_2(x)$,
$M_2(x)\uparrow\infty
$ with
$x$. The previous results, as well as Theorem
\ref{fr1}, Corollary
\ref{J}, suggest to study the behavior for 
$x$ large, of the averages 
 \begin{equation}\label{av}{1\over x}  \#\Big\{ m\le x:\inf_{ M_1(x)<N\le M_2(x)} \sup_{N\le
s_j^2\le
 Nf(N) }{ |
\o(m,j)- s_j^2 |\over s_j }\le z\Big\}.
\end{equation}
Here we set 
 \begin{equation} \label{sj}s_j^2: =\sum_{p\le j} {1\over p}-{1\over p^2}= \log\log j +\mathcal
O(1),
\end{equation}  
and the last relation comes from Mertens estimate.  For technical reasons (scale invariance properties of $W$ and   Kubilius model,
see next section), it turns up that it is  more convenient to replace the "$\log\log j$" term appeared before by $s_j^2$. 
 The resulting modifications  are thus neglectable in the statements. We show
that the  asymptotic order of the averages (\ref{av})  can be quantified  by  using Ornstein-Uhlenbeck process.    More precisely,  let
 \begin{equation}\label{in}I_N=I_N(f)=\big\{j: N\le s_j^2\le  Nf(N)\big\}.
\end{equation}
Let also $\mathcal N$ denote some increasing sequence of positive reals tending to infinity.  The theorem below allows to reduce the
study of the  averages (\ref{av})   to the one of     similar questions for  the Ornstein-Uhlenbeck process. Other
formulations may be easily extrapolated from the proof.  
\begin{theorem} \label{outom}  Assume that $M_2(x)= \mathcal O_\e (x^\e)$. Let $0<z''<z<z'$. As $x$ tends to infinity,
$$ \P\Big\{ \inf_{
M_1(x)<N\le M_2(x)\atop 
N\in \mathcal N}  \sup_{j\in I_N } { |W(s_j^2)  |\over s_j}\le z''\Big\} + o ( 1)\le \qq\qq\qq \qq \qq\qq\qq  $$
$$ {1\over x}  \#\Big\{ m\le x:\inf_{ M_1(x)<N\le M_2(x)\atop 
N\in \mathcal N} \sup_{j\in I_N } { | \o(m,j)- s_j^2 |\over s_j }\le z\Big\}  $$
$$\qq\qq\qq \qq\qq\qq \le \P\Big\{ \inf_{
M_1(x)<N\le M_2(x)\atop 
N\in \mathcal N}  \sup_{j\in I_N } { |W(s_j^2)  |\over s_j}\le z' \Big\} + o ( 1).
$$
 \end{theorem} 
By combining with   Corollary \ref{J}, we deduce for
instance 
 \begin{corollary} \label{oma} Assume that $M_2(x)= \mathcal O_\e (x^\e)$ and $   \log M_1(x)  = o\big( \log M_2(x)  \big)$. Let
$c<1/\l(z)$.   Then
$$\lim_{x\to \infty} {1\over x} 
\#\Big\{1\le  m\le x:\inf_{ \log M_1(x)<k\le \log M_2(x)}
\sup_{e^k\le s_j^2\le  e^kk^c} { |
\o(m,j)- s_j^2 |\over s_j }\le z\Big\}=1.  $$
\end{corollary}  
\begin{remark} \rm Let  
$d>1 $. There is no loss when restricting to   $x^{1/d}\le m\le x$ in the above ratios. But 
$  M_1(x)< e^k\le   M_2(x)$ imply $  M_1(m)< e^k\le   M_2(m^d)$. 
This allows to deduce from Corollary \ref{oma} a  result similar to those in \cite{FT} previously described, namely for
any
$d>1 $,
$$\lim_{x\to \infty} {1\over x} 
\#\Big\{1\le  m\le x:\inf_{ M_1(m)< e^k\le   M_2(m^d)}
\sup_{e^k\le s_j^2\le  e^kk^c} { |
\o(m,j)- s_j^2 |\over s_j }\le z\Big\}=1,  $$
  Taking for
instance  
$M_2 (x)=
\log x$,
$M_1 (x)= 
\log^{\e(x)} x$,
$\e(x)\to 0$, we get 
\begin{equation} \label{ldo1}   \lim_{x\to \infty} {1\over x}  
\#\Big\{ 1 \le m\le x:\inf_{   (\log  m)^{\e(m)}\le e^k\le   \log m}  \sup_{e^k\le s_j^2\le  e^kk^c }{ | \o(m,j)- s_j^2
|\over s_j }\le  z\Big\}=1.
\end{equation}  
And the relation $c\le 1/\l(z)$    asymptotically becomes $  \pi \sqrt c/2\le  z  $,   $z\to 0$. 
\smallskip\par
  \end{remark}
 \begin{remark} \label{ft}\rm  If   instead of condition   $   \log M_1(x)   = o\big( \log M_2(x)  \big)$, we have the weaker assumption $
\log M_1(x) = 
\rho\log M_2(x) 
 $,
$0<\rho <1$, then by operating similarly and using 0--1 law, we would also get for $\rho$ sufficiently small
 \begin{equation} \label{ldo1}   \lim_{x\to \infty} {1\over x}  
\#\Big\{ 1 \le m\le x:\inf_{   M_2(m)^\rho\le e^k\le   M_2(m)}  \sup_{e^k\le s_j^2\le  e^kk^c }{ | \o(m,j)- s_j^2
|\over s_j }\le  z\Big\} =1.
\end{equation}
However, we have no idea about a suitable precise value of $\rho$. 
  \end{remark}

 \medskip\par
 We will further establish a   delicate frequency result for the
truncated divisor function, which is in the spirit of Theorem
\ref{fr1}.
\begin{theorem} \label{lfro} Let   $ 0\le  M_1(x)<M_2 (x) $, $    M _1(x)\uparrow\infty$   such that 
  $M (x)= \mathcal O_\e (x^\e)$.   For any $z'>z>0$ and   $c<1/\l(z)$, there exists a constant $\k >0$ 
 such that,   $$ \lim_{x\to \infty} {1\over x} \#\bigg\{ m\le x:\! \inf_{ M_1(x)\le    n\le M_2 (x) }   {\!\#\Big\{  k\le n :
\displaystyle{\sup_{e^k\le s_j^2\le  e^kk^c}}   { | \o(m,j)- s_j^2
|\over s_j }\le  z' 
\Big\}\over   n^{1-c\l(z)}} 
 \ge \k \bigg\}$$$$= 1.
$$
  \end{theorem}

\section{Auxiliary results} We first list the needed probabilistic results. Next we briefly describe Kubilius model and extract from the
fundamental inequality a useful lemma.   The underlying small deviation  problem, namely the study for   small
$z
$ of 
$$\P\Big\{\sup_{0\le s\le t}|U(s)|<z\Big\},$$    can be yield to be 
  intimately linked to the Sturm-Liouville equation
\begin{equation} \label{sl/e} \psi''(x)-x\psi'(x)= -\l \psi(x), \qq \psi(-z)=\psi(z)=0.
\end{equation}
Let $\l_1\le \l_2\le \ldots$ and $\psi_1(x) , \psi_2(x),\ldots$ respectively denote the eigenvalues and normed eigenfunctions of Equation (\ref{sl/e}).  
Here $\l_i, \psi_j$ depend on $z$ and it is known that $\psi_1, \psi_2, \ldots$, form an orthonormal sequence with respect to the weight function
$e^{-x^2/2}$. According to Newell's result (see \cite{N} , see also (3.16) in \cite{C})
\begin{equation} \label{sl/n}\P\Big\{\sup_{0\le s\le t}|U(s)|<z\Big\}= {1\over (2\pi)^{1/2}}\sum_{k=1}^\infty
e^{-\l_k t}\Big(\int_{-z}^z\psi_k(x)e^{-x^2/2}\dd x\Big)^2.
\end{equation}
Let $\l(z)$ denote   the smallest eigenvalue  ($\l(z)=\l_1$). Then $\l(z)>0$ is a strictly decreasing continuous function of $z$ on
$(0,\infty)$. Further  
  \begin{equation} \label{sl/smei}\l(z)\sim {\pi^2\over 4z^2} \qq\qq {\rm  as}\   z\to 0 .
\end{equation}
   See Lemma 3.1 in \cite{C},  see also Lemma 2.2 for the following  result. 
  \begin{lemma} {\gsec (Cs\'aki's  estimate)} \label{lem1} For $z>0$, $t>0$ we have
 $$ {e^{-\l(z)t}\over (2\pi)^{1/2}}\Big(\int_{-z}^z\psi_1(x)e^{-x^2/2}\dd x\Big)^2\le \P\Big\{\sup_{0\le s\le t}|U(s)|<z\Big\}\le { e^{-\l(z)t}\over
1-e^{-t}}.
$$
  \end{lemma}
 
It follows  that for $z>0$, there exist positive constants $K_1(z), K_2(z)$ such that for all $t\ge 1$ 
\begin{equation}\label{est/ou/cs}K_1(z) e^{-\l(z)t}  \le \P\Big\{\sup_{0\le s\le
t}|U(s)|<z\Big\}\le K_1(z)  e^{-\l(z)t} . 
\end{equation} 
 \bigskip\par 
Now let   $\mathcal E_ s^t $  denote   the vector space generated $ U(u)$, $s\le u\le t$ and introduce the  
maximal correlation   coefficient      \begin{equation}\label{rc} \rho(\tau)= \sup_{ \xi\in \mathcal E_{-\infty}^{t}\atop
\eta\in\mathcal E_{t+\tau}^{\infty}}   \frac{ |\E(\xi-\E \xi)(\eta-\E \eta)|}{\big[|\E(\xi-\E \xi)^2\E(\eta-\E
\eta)^2|\big]^{1/2}} .
\end{equation}  
  By stationarity, this one does not depend on $t$. Stationary Gaussian processes such that $\rho(\tau)\to 0$ as $\tau\to \infty$  are called completely
regular. The spectral density of $U$ has the form $|\Gamma(\l)|^{-2}$ with $\Gamma(\l)= 1+i\l$, which is obviously  an entire function. Moreover, we
also have       
$\frac{\log |\Gamma(\l)|}{1+\l^2}\in L^1(\R)$. Further $\Gamma$ has $   i$ as unique imaginary zero. As $\Im\big(\frac{ 1}{\l -i}\big) =
\frac{ 1}{1+\l^2}$,
it follows that $\sup_{\l\in \R}\big|\Im\big(\frac{ 1}{\l -i}\big)\big|<\infty $.
Our next lemma is therefore just a direct consequence of  Theorem 6, section VI.6 in \cite{IR}.   
\begin{lemma}  \label{kr} The process $U$ is completely
regular, and further
$$\rho(\tau) =\mathcal O_\e(e^{- (1-\e)\tau}).$$
\end{lemma} 
   This result, which  is due to Kolmogorov and Rozanov (\cite{KR}, see  Theorem 1 and remarks at end of p.207),   will be crucial   in the proof of
Theorem
\ref{fr1}.  

\smallskip\par  Recall also the classical form of the Borel-Cantelli quantitative Lemma (\cite{P}, Theorem 3 or \cite{W}, Theorem 8.3.1).
    \begin{lemma}  \label{qbc}  Let 
$
\{A_k, k\ge 1\}$ be a sequence of events satisfying 
$$\P(A_k\cap A_\ell)\le   \P(A_k) \P(A_\ell) + \g_{\ell-k}
\P(A_\ell),\qq\quad (\forall \ell\ge k\ge 1)    $$
 where $ \g_i\ge 0 $ and   $\sum_{i=0}^\infty
\g_i<\infty$.
  Let $\psi_n=\sum_{k=1}^n \P(A_k)$ and assume that $\psi_n\to \infty$ with $n$. Then  for every $a>3/2$, 
 \begin{eqnarray*} \sum_{k=1}^n\chi_{A_k}&\buildrel{a.s.}\over{=}&\psi_n+ \mathcal O_a\big( \psi_n^{1/2}\big(
\log \psi_n)^a\big) 
 . 
\end{eqnarray*}  
\end{lemma}   
 
We finally need a suitable invariance principle for sums of independent random
variables. This one is due to Sakhanenko (see \cite{S}, Theorem 1). We give its most appropriate formulation for our
purpose. 
  Let $\{\xi_j, j\ge
1\}$ be independent centered random variables    with absolute second moments. Let $t_k= \sum_{j=1}^k \E \xi_j^2$,  
 $S_k=  \sum_{j=1}^k  \xi_j  $ and let $\{r_k, k\ge 1\}$ be some non-decreasing sequence of positive reals. Let $\a\ge 2$, $y>0$. Put 
successively,
\begin{eqnarray} \D_n&=&\sup_{k\le n}|S_k-W(t_k) |,\cr 
\D & = & \sup_{n\ge 1} {\D_n\over r_n},\cr 
\overline{\xi} &=& \sup_{j\ge 1} {|\xi_j|\over r_j},\cr
L_\a(y)&=& \sum_{j\ge 1} \E \min\Big\{{|\xi_j|^\a\over y^\a r_j^\a}, {|\xi_j|^2\over y^2 r_j^2}\Big\}.
\end{eqnarray}
 \begin{lemma} \label{sha1}There exists an absolute constant $C$ such that for any fixed $\a$, there exists a Brownian motion $W$ such that
for all
  $x>0$,
$$ \P\big\{\D  \ge C \a x\big\}\le L_\a(x) .  $$\end{lemma}
  
\bigskip\par Now we pass to the Kubilius model. Recall that  $p$   denotes some arbitrary prime number. Let $\{Y_p, p\ge 1\}$ be a
sequence of independent binomial random variables such that
$\P\{Y_p=1\}=1/p$ and
$\P\{Y_p=0\}=1-1/p$. We can view $Y_p$ as modelling whether or not an integer taken at random is divisible by $p$. Let 
$$ T_n=\sum_{p\le n}  Y_p, \qq\qq S_n=T_n -\E T_n.  $$
Then   $\E S_n^2= s_n^2  =\log\log n +\mathcal
 O(1)$ by (\ref {sj}).
The sequence $\{T_n, n\ge 1\}$ is known to asymptotically behave as the truncated prime divisor function
$$\o(m,t)=\#\{ p\le t : p|m\} ,  $$
  at least when $t$ is not too close to $m$.  More precisely, let
$$\o_r(m)=\big( \o(m,1), \ldots, \o(m,r)\big), $$
where $r$ is some integer with $2\le r\le x$, and put $u={\log x\over \log r}$. Then, given $c<1$ arbitrary, we have {\it uniformly} in
$x,r
$ and
$Q\subset \Z^r$,
\begin{equation} \label{km} {\#\{ m\le x: \o_r(m)\in Q\}\over x}=\P\big\{(T_1, \ldots, T_r)\in Q\big\}+ \mathcal O\big( x^{-c}+ 
e^{-u\log u}\big).
\end{equation}
See Lemmas 3.2, 3.5 in \cite{E} Chapter 3. See also \cite{Te}, Theorem 1 for a more precise result involving the Dickman function.
\begin{remark} \rm    There are natural   
  restrictions in the application of this estimate to {\it asymptotic} studies, due to the error term $e^{-u\log u}$.  To make it small,
it requires     if
$r=r(x)$ that
$  r(x)=\mathcal O _\e (x^\e)$   for all $\e>0$.   This amounts to truncate  the prime divisor function $\o (m)$  at level
$\mathcal O _\e (x^\e)$, which is satisfactory as long as $m\ll x$. However, these integers have a neglectable contribution on the size
of    the left-term of (\ref{km}).   Therefore the  model is mostly adapted to the analysis of the distribution of the small divisors
of an integer. See   \cite{E} p.122, see also     \cite{Te}  (Introduction) for a complete   and precise analysis of this point.  
  \end{remark}
  Estimate (\ref{km}) can be for instance used to   estimate the number of    
integers having no prime divisors in prescribed sets.
 Let $I =[p,q] $, $q\le x$; as $ \#\big\{m\le x:p|m  \Rightarrow  p\notin I\big\}  =  \#\big\{m\le x:\o(m,p)=\ldots =\o(m,q)\big\}$,
  it follows that
 \begin{eqnarray}{1\over x}\#\Big\{m\le x:p|m\ \Rightarrow  p\notin I\Big\}& =&  \prod_{p\in I}\big(1-{1\over p}\big)
 + \mathcal O\big( x^{-c}+  e^{-u\log u}\big). \end{eqnarray} 
 Choosing $I=[2,y]$, next $I=[ y,
x]$
allows to recover   known formula
 on  the   smallest
 or largest  prime divisors of $m$. 
\smallskip\par
Clearly, the approximation formula  (\ref{km}) can be used   to   transfer   properties from $(T_k)$ to
$\o$. Let indeed   $f$ be such that $f(N)= o_\rho(N^\rho)$.
 Recall that 
 $I_N =\big\{j: N\le
s_j^2\le
 Nf(N)\big\} $ and 
 let $\mathcal N$ be some fixed increasing sequence of reals. 
 Moreover, let $M_i :\N\to\R^+$   be non-decreasing  with $ \lim_{x\to\infty} M_i(x)=\infty$, $i=1,2$, and  such that 
 \begin{equation}\label{mi}1\le M_1(x)<M_2(x), \quad \qq  M_2(x)= \mathcal O_\e (x^\e). 
\end{equation}Let $r=r(x)\sim M_2(x)f(M_2(x))$, $r $ integer. Then 
$$ u=u(x)={\log x\over \log r(x)}\sim {\log x\over \log M_2(x)f(M_2(x))}\to \infty
 $$
with $x$.  Put 
 $$Q_{x}=\bigcup_{M_1(x)<N\le M_2(x)\atop 
N\in \mathcal N}\Big\{(\nu_1,\ldots,\nu_r)\in \Z^r:     \sup_{j\in I_N } { |\nu_j- s_j^2 |\over s_j}
 \le z\Big\}.$$
By  applying     (\ref{km}) with $Q=Q_{x}$, we get the useful comparison relation
\begin{lemma}\label{otos}   For any $z>0$,  as $x$ tends to infinity,
$${1\over x}  \#\Big\{ m\le x:\inf_{ M_1(x)<N\le M_2(x)\atop 
N\in \mathcal N} \sup_{j\in I_N } { | \o(m,j)- \log\log j |\over \sqrt{\log\log  j}}\le z\Big\}$$$$ =\P\Big\{ \inf_{
M_1(x)<N\le M_2(x)\atop 
N\in \mathcal N}  \sup_{j\in I_N } { |T_j - \log\log j |\over \sqrt{\log\log  j}}\le z\Big\}+ o ( 1).
$$
  \end{lemma}   
 
   
 \section{Proof of Theorem \ref{fr1}} 
     By stationarity  and by using
(\ref{est/ou/cs}), 
  \begin{equation}\label{est/ou/cs1}{K_1(z)\over f(e^k)^{ \l(z)} }
\le \P\big( A_k(f,z)\big)=\P\Big\{\sup_{0\le s\le  \log f(e^k) }|U(s)|<z\Big\}\le {K_2(z)\over  f(e^k)^{ \l(z)} } . 
\end{equation}
By summing up,
 \begin{equation}\label{est/ou/cs2} K_1(z) \sum_{k=1}^n  f(e^k)^{- \l(z)}\le \nu_n(f,z) \le K_2(z) \sum_{k=1}^n  f(e^k)^{- \l(z)} . 
\end{equation}  

If the series $\Sigma(f)= \sum_k f(e^k)^{- \l(z)}$ converges, by the first Borel-Cantelli lemma
$$\P\Big\{\sup_{k\le s\le k +\log f(e^k) }|U(s)|> z,  \qquad k  \ {\rm eventually}\Big\}=1. $$
Hence $f \in \mathcal U_z$. Now consider the case $\Sigma(f)=\infty$. We shall prove that $f \in \mathcal V_z$.  Let
$0<c_1<1/\l(z)<c_2$ and put 
$$f_1(t) =\log^{c_1} t, \qq f_2(t)=\log^{c_2} t.   $$
We may assume $f_1\le f\le f_2$. This is a standard device. Indeed, as $f_2\in \mathcal U_z$, we have the implication: $(f_1\vee f)\wedge
f_2\in \mathcal V_z\Rightarrow (f_1\vee f)\in \mathcal V_z\Rightarrow   f \in \mathcal V_z$. So it suffices to prove that $(f_1\vee
f)\wedge f_2\in \mathcal V_z$.  We now use the simplified notation
$ A_k(f,z)=A_k$ and notice that $ K'(z)  k^{-c_2\l(z)  }\le \P(A_k ) \le  K''(z)  k^{-c_1\l(z)  }$. By Lemma
\ref{kr}, 
 for  every  $ \ell>k$,  
$${\big|\P(A_k\cap A_\ell)- \P(A_k )\P(  A_\ell)\big|\over 
\sqrt{\P(A_k )(1-\P(A_k ))\P(  A_\ell)(1-\P(A_\ell ))}}\le C_1 e^{-C_2 ( \ell-k)} ,$$
$C_1, C_2$ being absolute constants. 
Hence
\begin{eqnarray*}\big|\P(A_k\cap A_\ell)- \P(A_k )\P(  A_\ell)\big|&\le &C_1 e^{-C_2 (\ell-k)} \sqrt{\P(A_k ) \P(  A_\ell) }\cr 
&\le &C_1 e^{-C_2 (\ell-k)} \P(  A_\ell)\Big({\P(A_k )\over  \P(  A_\ell)} \Big)^{1/2}
\cr &\le& C(z) \P(  A_\ell)  e^{-C_2 (\ell-k)}\Big( {\ell^{c_2 }\over  k^{c_1  }     
}\Big)^{\l(z)/ 2}  .\end{eqnarray*}

  But for some absolute constant  $C_3< C_2$ and  $C_4>0$ depending on $z$, we have
$$ e^{-C_2 (\ell-k)}\Big( {\ell^{c_2 }\over  k^{c_1  }     
}\Big)^{\l(z)/ 2} \le C_4 e^{-C_3 (\ell-k)}.  $$
 Indeed let
$ \ell= (H+1)k
$,
$H\ge 0$. This amounts to show that
$$    (H+1)      
  ^{c_2\l(z)/ 2}k^{(c_2-c_1)\l(z)/ 2} \le C_4 e^{ (C_2-C_3) Hk}.  $$
We use the following inequality. Let $\d,\b,\e$ be positive reals with $\d\ge \b$. Then there  exists $C$ depending on $\d,\b,\e$ only such that $H^\d
k^\b\le C e^{\e Hk}$ for all non-negative  reals $H,k$ with $k\ge 1$. Indeed, if $0\le H\le 1$, then 
 $H^\d k^\b\le (H  k)^\b\le C e^{\e Hk}$. And if $H>1$,  $H^\d k^\b\le (H  k)^\d\le C e^{\e Hk}$.

Applying this with $\d=c_2\l(z)/ 2$, $\b=(c_2-c_1)\l(z)/ 2$ yields
 $$     H       
  ^{c_2\l(z)/ 2}k^{(c_2-c_1)\l(z)/ 2} \le C e^{ \e  Hk},  $$
which implies our claim.
Thereby   
$$\big|\P(A_k\cap A_\ell)- \P(A_k )\P(  A_\ell)\big|\le C_4 e^{-C_3 (k-\ell)} \P(A_\ell).$$
 Lemma \ref{qbc} thus applies, and we deduce (for every $a>3/2$), 
 \begin{eqnarray*} \sum_{k=1}^n\chi_{A_k}&\buildrel{a.s.}\over{=}&\nu_n(f,z)+ \mathcal O_a  \Big(\nu_n(f,z)^{1/2} \log^a  \nu_n(f,z)\Big) 
 . 
\end{eqnarray*}
In particular   
$$\P\Big\{\sup_{k\le s\le k +\log f(e^k) }|U(s)|\ge z,  \qquad k  \ {\rm infinitely \ often}\Big\}=1. $$
Hence also $f\in \mathcal V_z$.  
 
Corollary \ref{J} follows easily. Indeed, let $0<c  \le 1/\l(z) $. By Theorem \ref{fr1}, $N_n(f_c,z)\uparrow\infty$ almost surely. And so
$\P\{  J(f_c ) \le z \}=1$. Now if $c> 1/\l(z) $, in view of estimate (\ref{est/ou/cs})  the series $\sum_{k=1}^\infty \P\{  A_k(f_c,z)\}$
converges. And by the first Borel-Cantelli lemma $\P\{ 
J(f_c ) > z \}=1$. 
 Corollary \ref{W}  is just a reformulation   of Corollary \ref{J} using the variable change $s= e^t$.

\section{ Proof of Theorem \ref{outom}}

   Now we can pass to the proof.  Let $\e, \eta  $ be positive reals. Let $\a$     sufficiently large so that $ \e\a> 1+\eta$.
Apply Lemma
\ref{sha1} to  
$S_n$   (here $\xi_p =Y_p-\E Y_p$).  Choose
$r_p=(\log\log  p)^{1+\eta\over
\a}$  and recall that 
  $\E |\xi_p|^\a  \sim  {1/p} $ for
  $p$  large.   Then 
$$\sum_{p }    {\E |\xi_p|^\a\over   r_p^\a}\le C \sum_{p}    {1\over  p (\log\log  p)^{1+\eta } }\le C \sum_{j}    {1\over  j\log j
(\log\log  j)^{1+\eta } } <\infty. $$
We have used the fact that if $p_j$ denotes the $j$-th prime number in the  increasing order, then $p_j \sim j \log j$. Now notice the
following simple estimate valid for all positive $y$,
$$L_\a(y)\le y^{-\a}\sum_{j\ge 1}    {\E |\xi_j|^\a\over   r_j^\a}. $$
 We deduce that $L_\a(y)\le C_\a y^{-\a}$. Recall that  $  \E S_n^2=s_n^2 $.  Therefore there exists a Brownian motion $W$
such that if 
\begin{equation} \label{M} {\bf \Upsilon}_\e   = \sup_{n\ge 1} {\sup_{j\le n}\big|S_j -W(s_j^2) \big|\over (\log\log  n)^{
\e}},\end{equation} 
 then $\E {\bf \Upsilon}_\e^{\b}<\infty$, $\b<\a$. We will just use the fact that $\E {\bf \Upsilon}_\e<\infty$. Let $z'>z$. By using
Lemma \ref{otos}, we have
\begin{eqnarray} \label{prelim}& &{1\over x}  \#\Big\{ m\le x:\inf_{ M_1(x)<N\le M_2(x)\atop 
N\in \mathcal N} \sup_{j\in I_N } { | \o(m,j)- s_j^2 |\over s_j }\le z\Big\}\cr 
&=& \P\Big\{ \inf_{
M_1(x)<N\le M_2(x)\atop 
N\in \mathcal N}  \sup_{j\in I_N } { |S_j   |\over s_j}\le z\Big\}+ o ( 1)\cr  &\le &\P\Big\{ \inf_{
M_1(x)<N\le M_2(x)\atop 
N\in \mathcal N}  \sup_{j\in I_N } { |W(s_j^2)  |\over s_j}\le z'\Big\}+ \P\{A\}+ o ( 1),
\end{eqnarray}
where we set
$$A=  \Big\{\inf_{
M_1(x)<N\le M_2(x)\atop 
N\in \mathcal N}  \sup_{j\in I_N } { |W(s_j^2)  |\over s_j}> z',  \inf_{
M_1(x)<N\le M_2(x)\atop 
N\in \mathcal N}  \sup_{j\in I_N } { |S_j  |\over s_j}\le z\Big\}. $$
We have 
$$\bigg|\inf_{
M_1(x)<N\le M_2(x)\atop 
N\in \mathcal N}  \sup_{j\in I_N } { |W(s_j^2)  |\over s_j}-  \inf_{
M_1(x)<N\le M_2(x)\atop 
N\in \mathcal N}  \sup_{j\in I_N } { |S_j  |\over s_j}\bigg|
$$ $$\le  \sup_{
M_1(x)<N\le M_2(x)\atop 
N\in \mathcal N}  \sup_{j\in I_N } { |S_j- W(s_j^2)     |\over s_j}\le  \sup_{
M_1(x)<N\le M_2(x) }  \sup_{j\le j^* } { |S_j- W(s_j^2)     |\over s_j} ,$$
where    $j^*$ denote  the largest indice such that $s_j^2\in  I_N$ of   $I_N $. 
Thus 
\begin{equation} \label{inter}\P\{A\}\le \P\Big\{ \sup_{
M_1(x)<N\le M_2(x) }  \sup_{j\le j^* } { |S_j- W(s_j^2)     |\over s_j}>z'-z\Big\}
.\end{equation}
\smallskip\par
 Let $\e'>\e$. Since $f(N)=o_\rho(N^\rho)$ by assumption and $s_j^2\sim \log\log j$ by (\ref{sj}), we have for all $N$ sufficiently
large, $N\ge N(\e, \e')$ say, 
$$\sup_{j\le j^*}\big|S_j -W(s_j^2) \big|\le {\bf \Upsilon}_\e\,(\log\log j^*)^{ \e} \le C{\bf
\Upsilon}_\e\,(N f(N))^{\e}\le C{\bf
\Upsilon}_\e\, N ^{\e'}.
  $$
  Then
\begin{equation}\label{apow}\sup_{j\le j^*}{ |S_j -W(s_j^2)  |\over s_j} \le CN^{- {1\over 2}}\sup_{j\le j^*}  |S_j -W(s_j^2)  | 
\le C{\bf \Upsilon}_\e\, N ^{ - {1\over 2}+\e'}
   .
\end{equation} 
Thereby for $N\ge N(\e, \e')$,
\begin{equation}\label{approxstow}\sup_{
M_1(x)<N\le M_2(x) }  \sup_{j\le j^* } { |S_j- W(s_j^2)     |\over s_j}\le  C {\bf \Upsilon}_\e\,M_1(x)^{-{1\over 2}+\e'} .
\end{equation} 
It follows that 
$$ \P\{A\}\le \P\Big\{ {\bf \Upsilon}_\e\, > C(z'-z)M_1(x)^{ {1\over 2}-\e'}\Big\}\le {C\over (z'-z)M_1(x)^{ {1\over 2}-\e'}}\, \E {\bf
\Upsilon}_\e.$$
Consequently $\P\{A\}= o(x)$, and we deduce from (\ref{prelim}) that
\begin{eqnarray} \label{prelim1}& &{1\over x}  \#\Big\{ m\le x:\inf_{ M_1(x)<N\le M_2(x)\atop 
N\in \mathcal N} \sup_{j\in I_N } { | \o(m,j)- s_j^2 |\over s_j }\le z\Big\}\cr  &\le &\P\Big\{ \inf_{
M_1(x)<N\le M_2(x)\atop 
N\in \mathcal N}  \sup_{j\in I_N } { |W(s_j^2)  |\over s_j}\le z'\Big\} + o ( 1).
\end{eqnarray}

 Now let $0<z''<z$.  As
$$  \P\Big\{\inf_{
M_1(x)<N\le M_2(x)\atop 
N\in \mathcal N}  \sup_{j\in I_N } { |W(s_j^2)  |\over s_j}\le  z'',  \sup_{
M_1(x)<N\le M_2(x)\atop 
N\in \mathcal N}  \sup_{j\in I_N } { |S_j-W(s_j^2)   |\over s_j}\le z-z''\Big\}. $$
$$\le \P\Big\{ \inf_{
M_1(x)<N\le M_2(x)\atop 
N\in \mathcal N}  \sup_{j\in I_N } { |S_j   |\over s_j}\le z\Big\},$$ 
we deduce from Lemma \ref{otos}
\begin{eqnarray} \label{prelim2}& &{1\over x}  \#\Big\{ m\le x:\inf_{ M_1(x)<N\le M_2(x)\atop 
N\in \mathcal N} \sup_{j\in I_N } { | \o(m,j)- s_j^2 |\over s_j }\le z\Big\}\cr 
&=& \P\Big\{ \inf_{
M_1(x)<N\le M_2(x)\atop 
N\in \mathcal N}  \sup_{j\in I_N } { |S_j   |\over s_j}\le z\Big\}+ o ( 1)\cr  &\ge &\P\Big\{ \inf_{
M_1(x)<N\le M_2(x)\atop 
N\in \mathcal N}  \sup_{j\in I_N } { |W(s_j^2)  |\over s_j}\le z''\Big\}-\P\{B\}+ o ( 1),
\end{eqnarray}
where we set
$$B=  \Big\{\!\!\inf_{
M_1(x)<N\le M_2(x)\atop 
N\in \mathcal N}  \sup_{j\in I_N } { |W(s_j^2)  |\over s_j}\le  z'',\!  \sup_{
M_1(x)<N\le M_2(x)\atop 
N\in \mathcal N}  \sup_{j\in I_N } { |S_j-W(s_j^2)   |\over s_j}> z-z''\Big\}. $$
By operating similarly, we also get 
\begin{eqnarray} \label{prelim3}& &{1\over x}  \#\Big\{ m\le x:\inf_{ M_1(x)<N\le M_2(x)\atop 
N\in \mathcal N} \sup_{j\in I_N } { | \o(m,j)- s_j^2 |\over s_j }\le z\Big\}\cr  &\ge &\P\Big\{ \inf_{
M_1(x)<N\le M_2(x)\atop 
N\in \mathcal N}  \sup_{j\in I_N } { |W(s_j^2)  |\over s_j}\le z''\Big\} + o ( 1).
\end{eqnarray}
 The proof is now complete.   

\begin{remark}\rm It follows from (\ref{apow}) that for all $0<\d<1/2$ and $d\ge 0$ 
   \begin{equation}\label{apow1}\E  \ \sup_{N}{1\over N^{ -\d }}  \sup_{j\in I_N}\Big|{  S_j -W(s_j^2)   \over s_j}\Big|^d 
<\infty .
\end{equation} 
   Consequently for any increasing unbounded sequence of reals $\mathcal N$,
  \begin{equation}\label{apow2} \liminf_{k\to \infty}\sup_{j\in
I_{N_k}}{|S_j|\over \sqrt{\log\log j}}\buildrel{a.s.}\over {=} 
\liminf_{k\to
\infty}\sup_{j\in I_{N_k}}{|W(s_j^2)|\over
\sqrt{\log\log j}}.
 \end{equation}
\end{remark}


  \section{Proof of Corollary \ref{oma}} Let $z''<z$. 
 Let $f=f_c$ with $c<1/\l(z'')$, $\mathcal N=\{e^k, k\ge 1\}$. Let also $0<\g<1$. Observe that 
\begin{eqnarray*} 
& &\P\Big\{ \inf_{
M_1(x)<N=e^k\le M_2(x) }  \sup_{j\in I_{N} } { |W(s_j^2)  |\over s_j}\le z''\Big\}
\cr &\ge & \P\Big\{ \inf_{
M_1(x)<N=e^k\le M_2(x) }  \sup_{e^k\le t\le  e^kk^c } { |W(t)  |\over \sqrt t}\le z''\Big\}
\cr &\ge &  \P\Big\{ \sum_{
\log M_1(x)<k\le \log M_2(x) } \chi\Big\{ \sup_{e^k\le t\le  e^kk^c } { |W(t)  |\over \sqrt t}\le z''\Big\}>0\Big\}\cr
&= &  \P\Big\{
\sum_{
\log M_1(x)<k\le \log M_2(x) } \chi \{ A_k(z'') \}>0\Big\}
 \cr&\ge & \P\Big\{ \sum_{
\log M_1(x)<k\le \log M_2(x) } \chi \{ A_k(z'') \}\ge \g\sum_{
\log M_1(x)<k\le \log M_2(x) } \P\{ A_k(z'') \} \Big\}.
\end{eqnarray*}
 Thus    
\begin{eqnarray}\label{comp}
& &\P\Big\{ \inf_{
M_1(x)<N=e^k\le M_2(x) }  \sup_{j\in I_{N} } { |W(s_j^2)  |\over s_j}\le z''\Big\}\qq\qq  \cr 
& & \qq\qq  =  \P\Big\{ {  N_{\log M_2(x)}(f_c ,z'')-N_{\log
M_1(x)}(f_c ,z'') \}\over \nu_{\log M_2(x)}(f_c ,z'')-\nu_{\log M_1(x)}(f_c ,z'')}\ge
\g 
\Big\}
.\end{eqnarray}

By Corollary \ref{J}, 
      
$$\lim_{n\to \infty}{N_n(f_c ,z'') \over \nu _n(f_c ,z'')}\buildrel{a.s.}\over{=}1 \quad {\rm and}  \quad K_1(z'')  \le  { \nu
_n(f_c,z)\over
n^{1-c\l(z'')}}
  \le  K_2(z'').$$

 By assumption, we have  $   \log M_1(x)  = o\big( \log M_2(x)  \big)$. Thus
$\nu_{\log M_1(x)}(f_c ,z'')= o\big(\nu_{\log M_2(x)}(f_c ,z'')\big)$. And it follows that 
  $$\lim_{n\to \infty} {  N_{\log M_2(x)}(f_c ,z'')-N_{\log M_1(x)}(f_c ,z'') \}\over \nu_{\log M_2(x)}(f_c ,z'')-\nu_{\log M_1(x)}(f_c
,z'')}\buildrel{a.s.}\over{=}1.$$
Consequently
$$\liminf_{x\to \infty}\P\Big\{ {  N_{\log M_2(x)}(f_c ,z'')-N_{\log M_1(x)}(f_c ,z'') \}\over \nu_{\log M_2(x)}(f_c ,z'')-\nu_{\log
M_1(x)}(f_c ,z'')}\ge
\g 
\Big\}=1.$$
By combining this with (\ref{comp}), we get
$$\liminf_{x\to \infty} \P\Big\{ \inf_{
M_1(x)<N=e^k\le M_2(x) }  \sup_{j\in I_{N} } { |W(s_j^2)  |\over s_j}\le z''\Big\}= 1$$
In view of Theorem \ref{outom}, this also implies  
 $$\lim_{x\to \infty} {1\over x}  \#\Big\{ m\le x:\inf_{ \log M_1(x)<k\le \log M_2(x)} \sup_{e^k\le s_j^2\le  e^kk^c} { | \o(m,j)- s_j^2
|\over s_j }\le z\Big\}=1.  $$
 The proof is now complete.

\section{Proof of Theorem \ref{lfro}}  The sets $A_k(c,z)$ being introduced before Theorem \ref{fr1}, we also define 
 \begin{eqnarray*}  
 B_k(c,z)&=&\Big\{ \sup_{j\in I_{N_k}}{ | \o(m,j)-
s_j |\over s_j}\le  z\Big\} 
 \cr C_k(c,z)&=&\Big\{ \sup_{j\in I_{N_k}}{ | S_j |\over
s_j}\le  z\Big\}
\cr D_k(c,z)&=&\Big\{ \sup_{j\in I_{N_k}} { |W(s_j^2)|\over s_j}\le 
z\Big\}   .\end{eqnarray*}
Fix $u>0$ and let $\eta>0$.
 By (\ref{apow1}), on a measurable set of full measure, we have for all $k$ large
enough,
 $ D_k(c,u ) \subseteq C_k(c,u+\eta ) $. 
Let $0<c< 1/\l(z) $. By Theorem \ref{fr1},
$$\lim_{n\to \infty} {\sum_{k=1}^n\chi_{A_k(c,z)}\over \sum_{k=1}^n\P(A_k(c,z))}\buildrel{a.s.}\over{=}1.$$
Let $0\le M_2(x)\uparrow \infty$   $ x$  and  such that 
 $M_2 (x)= \mathcal O_\e (x^\e)$. Obviously,
$$\lim_{x\to \infty} \inf_{M_1(x)\le n\le M_2(x)}{\sum_{k=1}^n\chi_{A_k(c,z)}\over \sum_{k=1}^n\P(A_k(c,z))}\buildrel{a.s.}\over{=}1.$$
Further
 $$ \k_1  \le  { \sum_{k=1}^n\P(A_k(c,z))\over n^{1-c\l(z)}}  \le  \k_2  ,$$
for some positive constants $\k_1, \k_2$. 
Let $z'>z$.  Since $A_k(c,z)\subseteq D_k(c,z') $, it follows that with probability one
\begin{eqnarray*}& &1 \buildrel{a.s.}\over{=} \lim_{x\to \infty} \inf_{M_1(x)\le n\le M_2(x)}{\sum_{k=1}^n\chi_{A_k(c,z)}\over
\sum_{k=1}^n\P(A_k(c,z))}\cr 
&\le&
\limsup_{x\to \infty} \inf_{M_1(x)\le n\le M_2(x)}{\sum_{k=1}^n\chi\{C_k(c,z')\}\over
\sum_{k=1}^n\P(A_k(c,z ))}   .
\end{eqnarray*} 
 Let $0<  \e<1$ and   put 
$$Q_{x}=\bigg\{(\nu_1, \ldots, \nu_r)\in \Z^r: \qq\qq\qq\qq\qq\qq\qq\qq $$
$$ \qq   \inf_{M_1(x) \le n\le M_2(x)}  \ {1\over \sum_{k=1}^n\P(C_k(c,z'))}\sum_{k=1}^n \chi\big\{\sup_{j\in I_{N_k}}{ |
\nu_j-s_j^2 |\over
s_j}\le  z'\big\} \le \e \bigg\}.$$
By applying (\ref{km}) with $Q=Q_{x}$, we get 
$$ {1\over x} \#\Big\{ m\le x: \inf_{M_1(x) \le n\le M(x)}  \ {1\over  
\sum_{k=1}^n\P(C_k(c,z'))}\sum_{k=1}^n \chi \{ B_k(c,z')\}  \le \e \Big\} $$ $$=\P\Big\{ \inf_{M_1(x) \le n\le M(x)}  \ {1\over  
\sum_{k=1}^n\P(C_k(c,z'))}\sum_{k=1}^n \chi \{ C_k(c,z')\}\le \e\Big\}+
o(1).
$$
Thus
$$ \limsup_{x\to \infty}{1\over x} \#\Big\{ m\le x: \inf_{M_1(x) \le n\le M(x)}  \ {1\over  
\sum_{k=1}^n\P(C_k(c,z'))}\sum_{k=1}^n \chi \{ B_k(c,z')\}  \le \e \Big\}$$
$$=\limsup_{x\to \infty}\P\Big\{ \inf_{M_1(x) \le n\le M(x)}  \ {1\over  
\sum_{k=1}^n\P(C_k(c,z'))}\sum_{k=1}^n \chi \{ C_k(c,z')\}\le \e\Big\} 
  =0 .
$$
This being true for all $0<\e<1 $, we  infer that
$$ \lim_{x\to \infty}{1\over x} \#\Big\{ m\le x: \inf_{M_1(x) \le n\le M(x)}  \ {1\over  
  \sum_{k=1}^n\P(C_k(c,z'))}\sum_{k=1}^n \chi \{ B_k(c,z')\}  \ge 1 \Big\} =1 .
$$
Finally, for some $ \k>0$ depending on $z$,
$$ \lim_{x\to \infty}{1\over x} \#\Big\{ m\le x: \inf_{M_1(x) \le n\le M(x)}  \ {1\over  
n^{1-c\l(z)}  }\sum_{k=1}^n \chi \{ B_k(c,z')\}  \ge \k  \Big\} =1 .
$$
   
\section{Proof of Theorem \ref{ind}}
Let $1/\a<\b<1/2$. Take $r_n=(\sum_{i=1}^j \E |X_i|^2)^{\b}=z_n^{2\b} $.  We notice   that
$$\sum_{j\ge 1} {\E |X_j|^\a\over r_j^{\a }}=\sum_{j\ge 1} {\E |X_j|^\a\over (\sum_{i=1}^j \E
|X_i|^2)^{\a\b}}\le C \sum_{j\ge 1} {\E |X_j|^\a\over (\sum_{i=1}^j \E |X_i|^\a)^{\a\b}}<\infty,$$
since $\a\b>1$. Thus
$$L_\a(y)\le y^{-\a}\sum_{j\ge 1} \E  {|X_j|^\a\over   r_j^\a} \le Cy^{-\a}.
$$
By Lemma \ref{sha1},   there exists   a Brownian motion $W$ such that if 
 $$  \Upsilon=\sup_{n}{1\over r_n}  \sup_{j\le
n}{|Z_j-W(z_j^2)| }    $$
then $\E  \Upsilon^{\a'}  <\infty$, $\a'<\a$. Now let $j_p^*= \max\{ j:  r_j\le 2^p\}$. As   
$$ \sup_{2^{p-1} <r_j\le 2^{p }
 }{|Z_j-W(z_j^2)|\over r_j }
 \le {2\over r_{j_p^*}}\sup_{ j\le j_p^*
 }{|Z_j-W(z_j^2)|  }  , $$
whenever  $\{j: 2^{p-1} <r_j\le 2^{p }
 \}\not=\emptyset$, it follows that 
$$    \sup_{ r_j\ge 1}
 {|Z_j-W(z_j^2)|\over r_j } \le 2\Upsilon. 
  $$
 Let
     $j(N)= \max(J_N)$.   
Hence 
\begin{eqnarray*} \Big|\sup_{j\in J_N}{|Z_j|\over z_j}-\sup_{j\in
J_N}{|W(z_j^2)|\over z_j}\Big|& \le & \sup_{j\in J_N}{|Z_j-|W(z_j^2)|\over z_j}=\sup_{j\in J_N}{|Z_j-|W(z_j^2)|\over
z_j^{1-2\b}z_j^{2\b}}
\cr &\le & \Big(\sup_{j\in J_N}  {1\over z_j^{  1-2\b }}\Big) \sup_{j\le j(N )}{|Z_j-|W(z_j^2)|\over  z_j^{2\b}} 
\cr &\le & \Big(\sup_{j\in J_N}  {1\over z_j^{  1-2\b }}\Big)\Upsilon  
\cr & \to &   0  ,
 \end{eqnarray*}
as $N\to \infty$ almost surely, since $\b<1/2$. By specifying this for $N=e^k$, we  therefore deduce
$$ \liminf_{k\to \infty}\sup_{e^k\le z_j^2\le e^kf(e^k)}{|Z_j|\over z_j}\buildrel{a.s.}\over{=}\liminf_{k\to \infty}\sup_{e^k\le
z_j^2\le e^kf(e^k)}{|W(z_j^2)|\over s_j}. 
$$
 almost surely.  This together with Corollary \ref{J} allows to conclude.

\section{Concluding Remarks}
Clearly, the    approximation formula (\ref{km}) applies to strongly additive  arithmetic   functions $f(n)=\sum_{p|n}
f(p)$, and associated truncated  functions. For additive  arithmetic   functions $f(n)=\sum_{p^\nu ||n }
f(p)$, the comparizon is made with the sums of independent random variables $\xi_p$ defined by $\P\{\xi_p= f(p^\nu)\} =(1-1/p)p^{-\nu}$,
$\nu=0,1,\ldots$. See \cite{E}, \cite{Te}. Special cases will be investigated elsewhere.

\medskip\par\noi
\noi{\sl Acknowledgments:} I am pleased to thank Endre Cs\'aki for the reference \cite{C}.
 
{
\end{document}